\documentclass[11pt]{amsart}
\usepackage{amsmath,amssymb}
\usepackage{setspace,xypic}

\newcommand\Ra{\mathbb R}

\newcommand\Za{\mathbb Z}

\newcommand\SL{\mathrm{SL}}
\newcommand\SO{\mathrm{SO}}
\newcommand\Sp{\mathrm{Sp}}

 \DeclareMathOperator{\Isom}{Isom}

\DeclareMathOperator{\CAT}{CAT}

\newcommand\ctz{\CAT(0)}
\newcommand\FHM{(FHM)}
\newcommand\FCP{(FCP)}
\newcommand\FH{(FH)}

\newtheorem{theorem}{Theorem}[section]
\newtheorem{proposition}[theorem]{Proposition}
\newtheorem{lemma}[theorem]{Lemma}
\newtheorem{corollary}[theorem]{Corollary}

\newtheorem{question}[theorem]{Question}
\newtheorem{conjecture}[theorem]{Conjecture}

\begin{document}

\title[Groups not acting]{Groups not acting on manifolds}
\author{David Fisher}
\address{Department of Mathematics, Indiana University,
Bloomington, IN 47401, USA}
\email{fisherdm@indiana.edu}
\thanks{The first named author partially supported by NSF grants DMS-0541917 and 0643546.}
\author{Lior Silberman}
\address{School of Mathematics, Institute for Advanced Study,
Princeton, NJ 08540, USA}
\email{lior@math.ubc.ca}
\thanks{The second named author partially supported by NSF grant DMS-0635607.}

\begin{abstract}
In this article we collect a series of observations that constrain
actions of many groups on compact manifolds.  In particular, we
show that ``generic" finitely generated groups have no smooth
volume preserving actions on compact manifolds while also producing
many finitely presented, torsion free groups with the same property.
\end{abstract}

\maketitle

\section{Introduction}
\label{section:introduction}

There are a number of interesting conjectures concerning actions of
large groups on manifolds, particularly conjectures of Gromov and
Zimmer on actions of higher rank lattices and Lie groups.  In this
context, Gromov conjectured that a random group should
not have any smooth actions on any compact
manifold.  In this paper we show that, in an appropriate
model of randomness, a random group has no smooth volume
preserving actions on compact manifolds.

We begin by defining the class of groups for which we can prove
this result.  That this class is in some sense ``generic" is
justified and discussed in Section \ref{section:generic}.  In that section we
also discuss some ``less generic" groups satisfying our
hypotheses.  While the notion of genericity we use necessarily produces
groups that are not finitely presented, we also provide many examples
of finitely presented groups satisfying our hypothesis.  In both cases,
we produce groups that are torsion free.
See \S \ref{subsection:questions} of this paper for further discussion
of both Gromov's conjecture and the meaning of ``generic" or ``random" group.

Let $\Gamma$ be a finitely generated group.  We say $\Gamma$ has
{\em no finite quotients} if there are no non-trivial
homomorphisms from $\Gamma$ to a finite group.  We say $\Gamma$
has {\em property $\FHM$} if any $\Gamma$ action on a complete
$\ctz$ Hilbert manifold has a fixed point.
By a {\em non-positively curved Hilbert manifold},
we mean a complete geodesic $\ctz$ metric space all of whose
tangent cones are (isometric to) Hilbert spaces. We remark that property $\FHM$
implies property $\FH$, the fixed point property on Hilbert spaces, which
(for locally compact groups) is equivalent to property $(T)$.

The main result of this paper is:

\begin{theorem}
\label{theorem:main} Let $\Gamma$ be a finitely generated group
with no finite quotients and property $\FHM$.  Then any volume
preserving action of $\Gamma$ on a compact manifold is trivial.
\end{theorem}

Since every finite group admits many actions on compact manifolds,
the assumption of no finite quotients is necessary. For weaker
statements on groups with property $\FHM$ but with finite quotients,
see Proposition \ref{proposition:fromzimmer} and Theorem \ref{theorem:finitequotients}.
We construct many torsion free groups satisfying the hypotheses of Theorem \ref{theorem:main}.
In Section \ref{section:generic},
we also discuss other classes of groups which satisfy the
conclusion of Proposition \ref{proposition:fromzimmer} without
having property $\FHM$. We include some observations concerning
groups with no actions by homeomorphisms on any compact manifold
in Section \ref{section:epilogue}. These last results depend
heavily on torsion elements.

The proof of Theorem \ref{theorem:main} involves three steps.
First, we observe that if a group $\Gamma$ has property $\FHM$ then any volume
preserving action on a compact manifold preserves
a measurable Riemannian metric.  Then we apply a theorem of Zimmer
\cite{Z2} to show that the invariant measurable metric and the fact
that $\Gamma$ has property $(T)$ imply that the action has discrete
spectrum, i.e. that the unitary representation of $\Gamma$ on
$L^2(M)$ decomposes as a sum of of finite dimensional subspaces.
The fact that $\Gamma$ has no finite quotients implies that it
has no non-trivial finite-dimensional representations.
It follows that the representation of $\Gamma$ on $L^2(M)$ is
trivial which immediately implies that the $\Gamma$ action on $M$ is trivial.

This article is motivated by the growing interest in many quarters
in the conjecture that random groups don't act on manifolds. No
one interested in the conjecture seemed to know the proof of
Theorem \ref{theorem:main} or it's application to ``generic"
finitely generated groups.

{\em Acknowledgements:}  The trick
of combining Zimmer's theorem from \cite{Z2} with no finite
quotients is first observed in  \cite{FM}, though in a slightly more
roundabout fashion.  Many thanks to Furman and Monod for interesting
conversations.

The torsion tricks used in section \ref{subsection:torsion} were
explained to the first author by Benson Farb in April of 2007. They
seem to have been observed by many people simultaneously and independently,
see e.g. \cite{BV,W}.
The application here to Kac-Moody groups appears to be new.

Thanks to Goulnara Arzhantseva and Ashot Minasyan
for pointing out the latter two methods of constructing groups without
finite quotients in section \ref{subsection:noquotients}.
Also thanks to Alain Valette and Yann Ollivier for useful remarks
on an earlier version of this paper and to Martin Bridson for sharing
the observation that other groups satisfy Corollary \ref{corollary:kacmoody}.

Conversations leading to this paper started in the workshop ``Geometrical
and Topological Rigidity" held at the Banff International Research Station
in July 2007.

\section{Proof of Theorem \ref{theorem:main}}
\label{section:proof}

We briefly recall the construction of the space of ``$L^2$
metrics" on a manifold $M$.  Given a volume form $\omega$ on $M$, we
can consider the space of all (smooth) Riemannian metrics on $M$ whose
associated volume form is $\omega$. This is the space of smooth sections
of a bundle $P{\rightarrow}M$.  The fiber of $P$ is $X=\SL(n,\Ra)/\SO(n)$.
The bundle $P$ is an associated bundle to the $SL(n,\Ra)$ sub-bundle of the
frame bundle of $M$ defined by $\omega$. The space $X$ carries a natural
$\SL(n,\Ra)$-invariant Riemannian metric of non-positive curvature; we denote its associated distance function by $d_X$.
This induces a natural notion of distance on the space of metrics, given by
$d(g_1,g_2)^2=\int_M d_X(g_1(m),g_2(m))^2 d\omega$.
The completion of the sections with respect to the metric $d$ will be denoted
$L^2(M,\omega,X)$;  it is commonly referred to as the {\em
space of $L^2$ metrics on $M$} and its elements will be called $L^2$ \emph{metrics} on $M$.
That this space is $\ctz$ follows easily from the fact that $X$ is $\ctz$.  For more
discussion of $X$ and its structure as a Hilbert
manifold, see e.g. \cite{FH}. It is easy to check that a volume preserving
$\Gamma$ action on $M$ defines an isometric $\Gamma$ action on
$L^2(M,\omega,X)$.

More generally, we can replace $X$ by any symmetric space $Y$ of
non-compact type and consider the same construction for any $Y$
bundle over $M$. In fact, the same construction applies if $(M,\omega)$
is just a standard finite measure space and does not depend on the
differentiable structure of $M$. The resulting space is called a {\em continuum
product}. One method for obtaining an isometric $\Gamma$ action on
$L^2(M,\omega,Y)$ is to have an $\omega$ preserving $\Gamma$ action on $M$
and a cocycle
$\alpha\colon\Gamma{\times}M{\rightarrow}\Isom(Y)$ satisfying
an integral bound (ensuring that the $\Gamma$ action preserves the space
$L^2(M,\omega,Y)$).  We call such actions {\em cocycle actions}.
Not all isometric $\Gamma$ actions arise in this
way.  This construction contains the case $Y=\Ra^n$, in which
case $L^2(M,\omega,Y)$ is a Hilbert space and there are many isometric
actions not arising
from cocycles over actions on $Y$.  This is essentially the only
way in which non-cocycle actions arise in this setting, see \cite{FH} for more discussion.

We say a group $\Gamma$ has property $\FCP$ if for any non-positively curved
symmetric space $Y$, any finite measure space
$(M,\omega)$ and any isometric
$\Gamma$ action on  $L^2(M,\omega,Y)$, the $\Gamma$ action has a
fixed point. Clearly property $\FHM$ implies property $\FCP$.
Since we do not assume our actions are cocycle actions,
any group with property $\FCP$ also has property $(FH)$.

Our argument establishes the following strengthening of
Theorem \ref{theorem:main}:

\begin{theorem}
\label{theorem:main'} Let $\Gamma$ be a finitely generated group
with property $\FCP$ and no finite quotients, then any volume
preserving $\Gamma$ action on a compact manifold is trivial.
\end{theorem}

\begin{lemma}
\label{lemma:invariantmetric} Let $\Gamma$ be a group with
property $\FCP$.  Then any volume preserving $\Gamma$ action on a
compact manifold $M$ preserves an $L^2$ metric.
\end{lemma}

This observation is immediate from the definitions and appears to be
well-known, but does not appear anywhere in the literature.

Combining Lemma \ref{lemma:invariantmetric} with a result of
Zimmer \cite[Theorem 1.7]{Z2}, we have:

\begin{proposition}
\label{proposition:fromzimmer}Let $\Gamma$ be a group with
property $\FCP$. Then any volume preserving $\Gamma$ action on a
compact manifold $M$ has discrete spectrum.
\end{proposition}

As mentioned in the introduction, {\em discrete spectrum} means
that the unitary representation of $\Gamma$ on $L^2(M,\omega)$
splits as an infinite direct sum of finite dimensional
representations.  In particular, Proposition
\ref{proposition:fromzimmer} implies that no group with property
$\FCP$ has a volume preserving weak mixing action on compact
manifolds.  As in \cite{Z2}, one can deduce from Proposition
\ref{proposition:fromzimmer} that the action is measurably
isometric, i.e. measurably conjugate to an action defined by
embedding $\Gamma$ in a compact group $K$.  Much stronger results
would follow if one could prove that $K$ was a Lie group.  To show
this, one can either show that enough of the $\Gamma$ invariant
subspaces of $L^2(M)$ are spanned by smooth functions or by
proving directly that the invariant metric is smooth (or even just
continuous).

Proving that $K$ is a Lie group is a well-known and difficult problem.
Here we do not need to establish this -- for our purposes it suffices to note
that it embeds in a product of compact Lie groups.
The proof of Theorem \ref{theorem:main'} is completed by the
following proposition.

\begin{proposition}
\label{proposition:nofiniteimages} Let $\Gamma$ be a group with no
finite images. Then any discrete spectrum action of $\Gamma$ is
trivial.
\end{proposition}

\begin{proof}
We have a $\Gamma$ action on $M$ whose action on $L^2(M,\omega)$
splits as a sum of finite dimensional representations $\pi_j$ on
finite dimensional spaces $V_j$.  Since finitely generated linear
groups are residually finite, each $\pi_j$ must have trivial
image.  Therefore the $\Gamma$ action on functions on $M$ is
trivial and so is the $\Gamma$ action on $M$.
\end{proof}

We remark briefly on one strengthening of our main results.
It is possible to have groups with only finitely many finite
quotients. For such a group $\Gamma$, there is always a maximal finite
quotient $F_{\Gamma}$.  Our methods also yield:

\begin{theorem}
\label{theorem:finitequotients}
Let $\Gamma$ be a group with property $\FCP$ and with finitely many
finite quotients.  Then any volume preserving $\Gamma$ action on a
compact manifold factors through $F_{\Gamma}$.
\end{theorem}

We remark that our results are stronger than the statements
of Theorem \ref{theorem:main'}, Proposition \ref{proposition:fromzimmer}
and Theorem \ref{theorem:finitequotients}.   Indeed, to obtain the conclusion
of those theorems, we only require a fixed point in any action on the space
$L^2(M,\omega, X)$ coming from a smooth action on $M$ or the even
weaker condition of a measurable invariant metric.  In particular,
the conclusion of Proposition \ref{proposition:fromzimmer} also holds for
all lattices in higher rank semisisimple algebraic groups over fields
of positive characteristic.  In this context the Zimmer-Margulis
approach to super-rigidity for cocycles produces a measurable
invariant metric and the groups are known to have property $(T)$,
so one can apply \cite[Theorem $1.2$]{Z2}. While it seems plausible that
these groups also have property $\FHM$, this does not seem to be known.

\section{Groups satisfying the assumptions of Theorem \ref{theorem:main}}
\label{section:generic}

In this section, we discuss methods of constructing groups satisfying the
hypotheses of Theorems \ref{theorem:main} and \ref{theorem:main'}.
In the first two subsections, we discuss groups with property $\FHM$
and property $\FCP$ respectively.  In the final subsection, we discuss
various methods which, starting with a hyperbolic group with
property $\FHM$ or $\FCP$, produce quotients of the given group with
no finite quotients.

\subsection{Groups with property $\FHM$.}
\label{subsection:fhm}  A criterion for property
$\FHM$ in terms of actions on simplicial complexes is given
explicitly in \cite{IN}.  This criterion builds on earlier work of
Wang \cite{W1,W2}.  A similar criterion is established in a
unpublished preprint of Schoen and Wang \cite{SW}.  Combined with
Zuk's work on random groups in the triangular model, this implies
that a random group in the triangular model at density more than
$1/3$ has property $\FHM$ with high probability.  As remarked in
\cite[I.3.g]{O} this then implies the same property for random
group in the density model with density more than $1/3$ with high
probability. More recently Naor-Silberman \cite{NS} have given
proofs that property $\FHM$ (and more) holds with high probability
for random groups in the graph model of \cite{Gr}. Also, Silberman
has given a simpler proof that property $\FHM$ holds in the density
model at density greater than $1/3$ \cite{S}. In the context of
\cite{NS}, we need much less than is used there.  Here we can get
by with adding relations corresponding to a single graph to a non-abelian free
group, rather than considering an infinite sequence of graphs. The
existence of such groups is fully justified by \cite{O2}.
We remark here that all the groups
mentioned in this paragraph are, with high probability, aspherical
and hyperbolic. This is important for constructing quotients of
these groups with no finite factors.  Here, a group is \emph{aspherical}
if it has an aspherical presentation. In particular, this implies the
group is torsion free.

Any cocompact group of isometries of a building of type $\tilde{A}_2$ can
be shown to have property $\FHM$ by the methods of \cite{IN}. These
groups therefore satisfy \ref{proposition:fromzimmer}.  It seems
quite likely that the same is true of cocompact groups of
isometries of irreducible higher rank buildings.  It seems
plausible that the same should be true for non-uniform lattices.
Since none of these groups is hyperbolic, we cannot use them to
build examples with no finite quotients.

\subsection{Groups with property $\FCP$.}
\label{subsection:fcp}
It follows from the main results of \cite{FH} that any quotient of
a lattice $\tilde \Gamma$ in $\Sp(1,n)$ by an infinite normal
subgroup has property $\FCP$. Again by standard constructions, one
can construct such a quotient $\Gamma$ which is torsion free and
hyperbolic, see \cite[Chapter II]{O} for discussion and
references.  It is not clear that one can construct the quotient
to be aspherical, so it is not clear that our first method for
producing groups with no finite quotients can be used for
these groups, see below.

\subsection{Groups with no finite quotients.\\}
\label{subsection:noquotients}

{\bf Method One:}
Given an aspherical hyperbolic group, there is a standard method
of using iterated random quotient to produce from it a group with
no fnite quotients. The finite stages of this process preserve the property
of being aspherical. The process is described on \cite[Section IV.k.]{O}
and is originally due to Gromov \cite{Gr0}.
This process can be applied to any aspherical hyperbolic group with property
$\FHM$ or $\FCP$, to obtain finitely generated, infinitely presented
groups with no finite quotients and with property $\FHM$ or
property $\FCP$.  This uses the fact that like property $\FH$,
properties $\FHM$ and $\FCP$ both obviously pass to quotients.
We remark here that this method of producing groups without
finite quotients depends only on having infinitely many relators
chosen at random.  The hard part of the construction is guaranteeing
that the resulting group is infinite.

It seems plausible that all assertions in the previous paragraph are true
for torsion free hyperbolic groups and not just aspherical ones. This is not
currently known and is technically a much more difficult question.  Because of
this difficulty, it is not clear that one apply this construction to the quotients of
lattices in $SP(1,n)$ in \S \ref{subsection:fcp}.

{\bf Method Two:} In \cite{Ol}, Ol'shanski gives a method of
producing infinite groups with no finite quotients from any
hyperbolic group.  It is not clear from that article that this
method can be used to produce torsion free groups though it may
be possible to achieve this using in addition arguments from
\cite{Ol1}.  The method does produce groups that are finitely presented.
This method applies to both random hyperbolic groups and to cocompact
lattices in $\Sp(1,n)$.

{\bf Method Three:}
In this method we use results from \cite{Os} and \cite{AMO} in a manner
inspired by \cite{ABJLMS}.

Let $F$ be any group with no finite quotients and at least three generators.
Then the free product $K=F*F$ is hyperbolic \emph{relative} to its two free
factors. Also, $K$ is finitely presented and torsion free if $F$ is.
Let $H$ be any relatively hyperbolic group.
Then by \cite[Theorem 1.4]{AMO} and
\cite[Theorem 2.4]{Os} there exists an infinite relatively hyperbolic group $\Gamma$ that is a
common quotient of $K$ and $H$ such that $\Gamma$ is finitely presented and torsion free
if $K$ and $H$. The peripheral subgroups of $G$ are exactly images of the
peripheral subgroups of $K$ and $H$.

To apply this result in our context it suffices to find groups
$F$ which are finitely presented, torsion free and have no finite quotients.
One can use, e.g. the finitely
presented torsion free simple groups constructed
by Burger and Mozes \cite{BM} or the four generated
Higman group with the presentation

$$\langle a,b,c,d | a^b =a^2,  b^c= b^2, c^d= c^2, d^a= d^2 \rangle.$$

This technique can be applied to produce groups with no finite quotients
from random hyperbolic groups and from both uniform and non-uniform lattices
in $\Sp(1,n)$.

When applying either of the last two methods to lattices in $\Sp(1,n)$, one can
apply them directly to the lattice.  Since the method produces infinite groups
with no finite quotients, it is clearly producing quotients of the lattice
by infinite normal subgroups of infinite index.

\section{Other groups not acting and some questions}
\label{section:epilogue}

We close with two collections of remarks.  In the first subsection,
we discuss other reasons why torsion impedes non-trivial group actions,
mostly to explain our emphasis on producing torsion free groups.  In
the second subsection, we ask some questions motivated by this work
and Gromov's conjecture.

\subsection{Torsion tricks}
\label{subsection:torsion}

This section points out some classes of groups with no actions by
homeomorphisms on compact manifolds.  The main point is the
following fact.  We learned it from Farb and Whyte, but it
appears to be known independently by many people see e.g. \cite{BV,W}.

\begin{lemma}
\label{lemma:criterion} Let $\Gamma$ be a simple group that
contains a copy of $(\Za / p\Za)^{\infty}$.  Then $\Gamma$ has no
non-trivial actions by homeomorphisms on any compact manifold.
\end{lemma}

The lemma is simply the fact, see e.g. \cite{MS}, that for any
compact manifold $M$, there is a number $k=k(M)$ such that a
faithful action of $(\Za / p\Za)^{n}$ by homeomorphisms on $M$
implies $n \leq k$.  Farb and Whyte observed that by results in \cite{Sch}
one can construct a two generated simple group
$\Gamma$ containing the lamplighter group $\Za \wr \Za/ p \Za$.  The lemma then
implies that this group $\Gamma$ has no non-trivial actions by
homeomorphisms on any compact manifold.  One can also give a proof
of Shmuel Weinberger's observation \cite{Z3} that $\SL(\infty,
\Za)$ does not act on a compact manifold using the same
observations, Margulis' normal subgroups theorem and the
congruence subgroup property for $\SL(n,\Za)$.

A more constructive method for finding finitely presented simple
groups containing $(\Za / p\Za)^{\infty}$ is the theory of
Kac-moody groups.  Lemma \ref{lemma:criterion}, the main theorem of \cite{CR} and
\cite[Proof of Theorem 4.6]{R} imply that:

\begin{corollary}
\label{corollary:kacmoody} Let $\Gamma$ be a split or almost split
Kac-Moody group over the finite field $\Za/p\Za$ with an infinite,
irreducible, non-affine Weyl group.  If $p$ is large enough, then
any $\Gamma$ action by homeomorphisms on a compact manifold is
trivial.
\end{corollary}

We remark here that the $\Gamma$ in the corollary is finitely
presented.  Bridson has independently observed that Thompson's
group $V$ and the variants constructed by Higman satisfy a similar
conclusion for the same reasons.  Kac-Moody groups are
of particular interest in this context because they are lattices
in locally compact groups.  Historically the motivation for
studying questions concerning groups (not) acting on compact
manifolds derives from Zimmer's conjectures concerning actions of
lattices in Lie groups and also algebraic groups over other local
fields.

In a recent paper \cite{ABJLMS}, similar and more elaborate uses of torsion have
led to classes of infinite groups with no non-trivial actions on certain classes
of non-compact manifolds and also more general spaces.

\subsection{Further questions}
\label{subsection:questions}

The results in this paper raise an obvious sequence of questions.
Are there finitely generated or even finitely presented, torsion
free groups with no non-trivial actions on compact manifolds? One
can ask the question either for actions by homeomorphisms or for
actions by diffeomorphisms.  It seems plausible that most groups
will have no smooth actions.  Unless the manifold acted on is the
circle, there are no known obstructions to any torsion free group acting by
homoeomorphisms on any compact manifold.  While producing a few sporadic examples would already be quite
interesting, there is significantly more interest in showing that having
no non-trivial actions is typical (``generic'') in a model for random groups.
However, it is important to note that groups chosen according to such a model
are not ``generic'' in the way that ``random regular graphs'' are --
there is currently no good notion of a ``typical'' group.

More concretely, to reduce Gromov's conjecture to Theorem
\ref{theorem:main}, one wants to prove:

\begin{conjecture}
\label{conjecture:volume} In an appropriate model for random groups,
any action of a``generic'' group $\Gamma$ on a compact manifold $M$
preserves a smooth volume form.
\end{conjecture}

The conjecture is interesting in the context of groups satisfying
Theorems \ref{theorem:main}, \ref{theorem:main'},
\ref{theorem:finitequotients} or even just
Proposition \ref{proposition:fromzimmer}.
It seems reasonable to try to approach
Conjecture \ref{conjecture:volume} as a fixed point problem.  It
is worth noting that the fixed point property used must be strictly
stronger than property $(T)$ as many linear $(T)$ groups admit
actions on manifolds that do not preserve a volume form.  A version
of this conjecture was asked as a question by Nigel Higson at the
July 2007 Banff workshop.

Another reasonable and fairly well known question also motivated
by the work in this paper and Gromov's conjecture.

\begin{question}
\label{question:finitegroups}
Is there a number $\frac{1}{2}>d>0$ for which a random group with density
at least $d$ has no finite quotients with positive probability?
With probability tending to $1$?
\end{question}

In particular, a positive answer to this question would produce many
hyperbolic groups which are not residually finite.

\bigskip

\end{document}